\documentclass[twoside,11pt]{article}

\usepackage{graphicx, subfigure}
\usepackage[top=1in, bottom=1in, left=1in, right=1in]{geometry}
\usepackage[sort&compress,numbers]{natbib} 
\usepackage{amsfonts, amsmath, amssymb, amsthm, constants}
\usepackage{MnSymbol}
\usepackage[mathscr]{eucal}
\usepackage{enumitem}
\setlist{noitemsep}
\usepackage{booktabs}
\usepackage{microtype}

\usepackage{titlesec}
\titleformat{\section}[block]{\large\bfseries}{\thesection.}{1em}{}
\titleformat{\subsection}[runin]{\bfseries}{\thesubsection.}{1em}{}

\usepackage{color}
\definecolor{darkred}{RGB}{100,0,0}
\definecolor{darkgreen}{RGB}{0,100,0}
\definecolor{darkblue}{RGB}{0,0,150}

\usepackage{hyperref}
\hypersetup{colorlinks=true, linkcolor=darkred, citecolor=darkgreen, urlcolor=darkblue}
\usepackage{url}

\urldef\eacurl\url{http://www.math.ucsd.edu/~eariasca/}

\def\ppv{{\sf PPV}}
\def\fdr{{\sf FDR}}
\def\fdp{{\sf FDP}}
\def\mfdr{{\sf mFDR}}
\def\fwer{{\sf FWER}}
\def\H{\mathscr{H}}
\def\taubh{\tau_\gamma}

\newtheorem{prp}{Proposition}

\theoremstyle{remark}

\newtheorem{rem}{Remark}

\def\beq{\begin{equation}} % \setcounter{equation}{1}}
\def\eeq{\end{equation}}
\def\beqn{\begin{eqnarray*}}
\def\eeqn{\end{eqnarray*}}
\def\Bitem{\begin{itemize}\setlength{\itemsep}{.2in}}
\def\bitem{\begin{itemize}\setlength{\itemsep}{.05in}}
\def\eitem{\end{itemize}}
\def\Benum{\begin{enumerate}\setlength{\itemsep}{.2in}}
\def\benum{\begin{enumerate}\setlength{\itemsep}{.05in}}
\def\eenum{\end{enumerate}}
\def\bmult{\begin{multline*}}
\def\emult{\end{multline*}}
\def\bcenter{\begin{center}}
\def\ecenter{\end{center}}
\def\bframe{\begin{frame}}
\def\eframe{\end{frame}}

\newcommand{\prpref}[1]{Proposition~\ref{prp:#1}}

\newcommand{\tabref}[1]{Table~\ref{tab:#1}}

\def\cN{\mathcal{N}}

\def\cR{\mathcal{R}}

\newcommand{\E}{\operatorname{\mathbb{E}}}
\renewcommand{\P}{\operatorname{\mathbb{P}}}

\def\Unif{\text{Unif}}

\def\implies{\ \Rightarrow \ }

\def\comp{\mathsf{c}}
\def\1{\mathbbm{1}}

\definecolor{purple}{rgb}{0.4,.1,.9}

\begin{document}

%\twocolumn

\thispagestyle{empty}

\title{Noisy Hypotheses in the Age of Discovery Science}
\author{Ery Arias-Castro\thanks{University of California, San Diego --- \eacurl}}
\date{}
\maketitle

\vspace{-0.3in}
\bcenter\small
\begin{tabular}{p{0.8\textwidth}}
%\toprule
%{\em Abstract.}
We draw attention to one specific issue raised by~\citet{ioannidis2005most}, that of very many hypotheses being tested in a given field of investigation.  To better isolate the problem that arises in this (massive) multiple testing scenario, we consider a utopian setting where the hypotheses are tested with no additional bias.  We show that, as the number of hypotheses being tested becomes much larger than the discoveries to be made, it becomes impossible to reliably identify true discoveries.  This phenomenon, well-known to statisticians working in the field of multiple testing, puts in jeopardy any naive pursuit in (pure) discovery science. \\  
%\bottomrule
\end{tabular}
\ecenter

%\tableofcontents

\section{Introduction} \label{sec:intro}

More than ten years ago now,~\citet{ioannidis2005most} declared that ``most published research findings are false".  If anything, his claim has been corroborated by the empirically observed lack of replicability in the sciences at large.  He provided the following argument to support his claim.  Consider an ideal setting where, in a certain field of scientific investigation, a number of hypotheses are tested.  All tests are assumed to be calibrated at level $\alpha$ and have power $1-\beta$ for the underlying effect.  Let $\pi_1$ denote the proportion of false null hypotheses among all the null hypotheses being tested, that is, the proportion of real discoveries to be made.  ~\citet{ioannidis2005most} defines the positive predictive value (PPV) as the probability that a given rejection is correct, which he derives under the present setting and obtains 
\beq\label{ppv}
\ppv{} = \frac{(1-\beta) \pi_1}{(1-\beta) \pi_1 + \alpha (1-\pi_1)}.
\eeq
The level $\alpha$ is typically set following standard practice ($\alpha = 0.05$ still being popular in some fields).  The power $1-\beta$ is increasing in $\alpha$, in the effect size, and in the sample size.  The proportion of false nulls $\pi_1$ is a consequence of how the hypotheses to be tested are chosen or formulated, and thus arises `organically' from the community of researchers interested in the field of expertise under consideration.
From this formula alone~\citet{ioannidis2005most} derives three corollaries: the first two are about the fact the PPV is increasing in the power and the third one on the fact that the PPV is increasing in the proportion of false nulls ($\pi_1$).
He encapsulates the impact of the latter in the following statement:
\beq\begin{array}{c} \label{ioannidis}
\text{``The greater the number and the lesser the selection of tested relationships} \\
\text{in a scientific field, the less likely the research findings are to be true."}
\end{array}\eeq
This is very easily seen from the fact that
\beq\label{ppv-ub}
\ppv \ge 1/2 \implies \frac{\pi_1}{1-\pi_1} \ge \alpha.
\eeq
To ensure that at least half of the claimed discoveries are correct, we see that $\alpha$ needs to be of order at most $\pi_1$ (and this is assuming a power of 1).  Thus, if $\pi_1$ is (very) small, $\alpha$ needs to be set (very) small, and this presents a problem as it lowers the power.  (Remember that the power of a test is decreasing in the level.)

The parameter $\pi_1$ could be called the ``probing coefficient".  A large $\pi_1$ (close to 1) is indicative of a field where hypotheses are generated with care, as might have been the case up to a few decades ago, where data were much more scarce and much more costly to collect and process, and also where the number of researchers (and perhaps the pressure to publish?) was much more limited.  
A small $\pi_1$ is evidence that hypotheses are generated (and tested) with abandon, which might be caricatural of the various `omics', where massive testing is commonplace.  
%This does not come without consequences and leads to unsettling conclusions.

%Suppose we want to control the PPV at $1-\gamma$ for some $\gamma \in (0,1)$.  Straightforward algebra yields
%\beq
%\ppv{} \ge 1-\gamma \iff \pi_1 \ge \frac{\alpha(1-\gamma)}{\alpha(1-\gamma) + (1-\beta)\gamma}.
%\eeq
%Assuming $\gamma = 0.10$, as well as $\alpha = 0.05$ and $1-\beta = 0.80$ (a very optimistic figure), the requirement is approximately $\pi_1 \ge 0.36$.  
%This is routinely violated in a single omics experiment.  

\citet{ioannidis2005most} derives the PPV in two other settings: when some form of bias is present and when multiple labs are pursuing the same question independently in parallel.  The corresponding formulae lead to three additional corollaries: two are about the fact that PPV is decreasing in the amount of bias and one about the fact that PPV is decreasing in the number of labs performing the same test in parallel.\footnote{In~\cite{ioannidis2005most} it is assumed that only the first lab reaching significance gets to publish its result.}
It is generally recognized that bias is an important factor in statistical analyses, due in part to a general preference for positive results (i.e, rejections).

Our main concern here is in exploring the effect that $\pi_1$ has on the PPV and its impact on the foundations of so-called `discovery science'.  
In the process, we place the PPV in the context of multiple testing, where it truly belongs, and discuss the issue of distinguishing false from true null hypotheses in this framework.

In a cartoonish simplification, the discovery science paradigm posits the availability of large quantities of data and the scientist --- now sometimes referred to as a `data scientist' --- lets the data speak for themselves to find meaningful effects.   
Our message in this paper is that, taken to an extreme, this paradigm is bound to fail if hypotheses are tested with abandon.  Taking a multiple testing perspective and drawing from the corresponding literature, we offer some arguments in support of the following statement:
\begin{quote}
{\em In a field of pure discovery science where very many hypotheses are tested, there is a risk that false nulls are indistinguishable from true nulls.}
\end{quote}

We argue our point by showing that, as the number of hypotheses increases, the threshold for significance needs to become ever closer to 0, making rejection of both true and false null hypotheses ever harder. 
This is rather intuitive.  Indeed, assuming the null P-values are uniformly distributed in $[0,1]$ under their respective nulls, and if there are a large number $m_0$ of true null hypotheses, the $k$-th largest will be close to $k/m_0$.  Qualitatively speaking, the smallest P-values from the true null hypotheses are ever smaller as their number increases, and the P-values from the false null hypotheses are `competing' with that.  We provide some elements of theory that go beyond this simple argument.

\section{A utopian model for pure discovery science}
To emphasize our point we consider what many would deem a utopian setting where, in a given field of investigation, a number of studies are performed.  We assume that each study is preregistered, conducted exactly as planned, and published without delay.  We also suppose that each study is based on the same sample size, that each study results in a single bonafide P-value,\footnote{If $p$ denotes the P-value, this typically means that $\P(p \le \alpha) \le \alpha$ for all $\alpha \in [0,1]$ under its null hypothesis.  Here we make the simplifying assumption that $p$ is uniformly distributed in $[0,1]$ under its null hypothesis.} that the sought-after effect size is the same, and that the sample size is large enough to result in sufficient power.
Moreover, in line with a pure discovery science perspective, we further assume that all studies are, a priori, equally important.  
\begin{table}[h]
\centering\small
\caption{The setting we consider versus the real state-of-affairs.}
\label{tab:setting}
\bigskip
\setlength{\tabcolsep}{0.22in}
\begin{tabular}{p{0.42\textwidth} p{0.42\textwidth}}
\toprule
%\multicolumn{2}{c}{Item} \\
%\cmidrule(r){1-2}
{\bf Utopia (our setting)} & {\bf Reality (state-of-affairs)} \\ \midrule
Each study is preregistered and conducted as planned & The preregistering of studies is not yet common practice~\cite{open2012open,dickersin2003registering} \\ \midrule
Each study is published & Publication bias can be important~\cite{dwan2008systematic} \\ \midrule
Each study results in a single test & Multiple testing within a single study is common (and sometimes massive)~\cite{dudoit2007multiple,rice2008methods} \\ \midrule
Each study results in a bonafide P-value & P-values are often approximate or based on unverifiable assumptions, and can be severely biased by informal explorations of the data (data snooping) or flawed by misuse of statistics~\cite{simmons2011false,gelman2013garden} \\ \midrule
Each study is well-powered & Studies are routinely underpowered~\citep{button2013power}\\ \midrule
Each study tests a different hypothesis & Multiple labs tackle the same question, possibly in different variants and using different tools \\ \midrule
Field-wide multiple testing is performed to identify the most significant studies & Published P-values are rarely considered in a larger context (except, possibly, in systematic reviews, but years later) \\
\bottomrule
\end{tabular}
\end{table}
We also postulate a fair distribution of resources.
Indeed, we imagine the whole field publishing all the results in a single repository and, at the end of each year (say), performing a meta-analysis of all the results published that year.  Since in the present setting the P-values carry equal weight, the field as a whole allocates the `follow up' budget equally to further investigate the studies with the most significant P-values.  
All this is summarized in \tabref{setting}, where we also contrast this setting with the reality of many fields of investigation.
Of course, the field as a whole would also allocate part of its budget to the generation of new data, new hypotheses, etc, but this will not be of interest here.

In this stylized setting, the sole issue one needs to contend with is the number of hypothesis tests being performed.  We implicitly assume this number is large given that we intend to model a field that relies heavily on a data-driven discovery paradigm.  Therefore, what remains to be done is to solve the multiple testing problem posed by the very many P-values that are generated.  
We assume there are $m$ null hypotheses being tested in the period of interest (we said one year above).  We denote by $\H_1, \dots, \H_m$ these hypotheses and $p_1, \dots, p_m$ the resulting P-values.  (We leave implicit the test statistics that generated these P-values.)  We let $m_0$ and $m_1$ denote the number of true and false nulls, respectively, and $\pi_0 = m_0/m$ and $\pi_1 = m_1/m$ their proportions.

\section{Two interpretations for the PPV}
We start with the PPV and derive it in two ways.  
The first derivation relies on Bayes' formula and is essentially the derivation of~\citet{ioannidis2005most}, although here we offer a frequentist perspective.  
The second derivation puts the PPV in the context of multiple testing, where we believe it belongs.  
%This may not have been known to~\citet{ioannidis2005most}.

\subsection{The derivation of the PPV from Bayes' formula}
\citet{ioannidis2005most} derives the formula for the PPV, presented in \eqref{ppv}, using Bayes' formula.  
We provide some details here not provided in the original article.  We note that Ioannidis takes a sort of Bayesian perspective when deriving the PPV, as he (implicitly) assumes that the status of each hypothesis (true or false) is random.  We offer a frequentist perspective instead.

Ioannidis postulates that any P-value not exceeding $\alpha$ is deemed significant (and the corresponding null hypothesis is therefore rejected).  
Suppose that we pick a P-value at random among the P-values that were found significant.  We ask ourselves the question: what is the probability that the corresponding null hypothesis is indeed false?  

We assume that each P-value is uniformly distributed in $[0,1]$ under its null hypothesis, so that
\beq\label{level}
\P(p_i \le \alpha) = \alpha, \quad \text{when $\H_i$ is true.}
\eeq
We also assume that all tests have the same power $1-\beta$, so that 
\beq\label{power}
\P(p_i \le \alpha) = 1-\beta, \quad \text{when $\H_i$ is false.}
\eeq
Let $I$ be a uniformly distributed index in $\{1, \dots, m\}$.
The question we are asking requires the computation of $\P(\H_I^\comp \mid p_I \le \alpha)$ where for each $i$, $\H_i^\comp$ indicates that $\H_i$ is false.  Using Bayes' formula, we derive
\beq
\P(\H_I^\comp \mid p_I \le \alpha) = \frac{\P(p_I \le \alpha \mid \H_I^\comp) \P(\H_I^\comp)}{\P(p_I \le \alpha \mid \H_I^\comp) \P(\H_I^\comp) + \P(p_I \le \alpha \mid \H_I) \P(\H_I)}.
\eeq
Since the proportions of true and false nulls are $\pi_0$ and $\pi_1$, 
\beq
\P(\H_I) = \pi_0 = 1 - \pi_1, \quad  \P(\H_I^\comp) = \pi_1.
\eeq
Also, by \eqref{level},
\beq
\P(p_I \le \alpha \mid \H_I) = \alpha,
\eeq
and by \eqref{power}, 
\beq
\P(p_I \le \alpha \mid \H_I^\comp) = 1 - \beta.
\eeq
Plugging this into the large fraction above results in the expression \eqref{ppv}.  Thus, we have established the following frequentist interpretation of the PPV.

\begin{prp}\label{prp:ppv}
Assume we are performing a number of tests at level $\alpha$, all having power $1-\beta$.  The PPV is the probability that a randomly chosen hypothesis among those rejected is false.
\end{prp}

We note that the independence of the P-values was not used.

\begin{rem}
Although the setting may appear to be stylized, a similar result holds much more generally.  Indeed, suppose that test $i$ is performed at level $\alpha_i$ and has power $1-\beta_i$.  The reader is invited to verify that \prpref{ppv} still holds in this setting, with \eqref{ppv} unchanged but with $\alpha$ defined as the average level and $1-\beta$ defined as the average power, meaning 
\beq
\alpha = \frac1m \sum_{i=1}^m \alpha_i, \quad
1-\beta = \frac1m \sum_{i=1}^m (1-\beta_i).
\eeq
\end{rem}

\subsection{The PPV as one minus the marginal false discovery rate (mFDR)}
We present another interpretation of the PPV that puts it in the larger context of multiple testing.  It turns out that the PPV can be interpreted as one minus the marginal false discovery rate (mFDR)~\cite{storey2007optimal} of the naive multiple test that rejects every null hypothesis with a P-value below $\alpha$ --- see \eqref{multiple-alpha} below.  We provide some background for the reader not already familiar with multiple testing.

In the particular setting considered here, where the P-values are assumed independent, we may without loss of generality define a multiple test as a function of the P-values $p_1,\dots, p_m$ that returns an index $\cR \subset \{1, \dots, m\}$ indicating the hypotheses to be rejected.  In formula, the procedure rejects $\H_i$ if $i \in \cR$.  
%($\cR$ is meant to indicate ``rejection".)  
Applying a given multiple test results in the situation described in \tabref{multiple}. 

\begin{table}
\centering
\caption{Some notation commonly used in multiple testing.  All the quantities, except for $m, m_0, m_1$ are random and depend on the (implicit) multiple testing procedure in use.  In this notation, $\pi_0 = m_0/m$ is the proportion of true nulls and $\pi_1 = m_1/m$ is the proportion of false nulls.}
\label{tab:multiple}
\medskip
\def\arraystretch{1.2}
\begin{tabular}{r|cc|c}
& True Null & False Null & Total \\ \hline
Rejection & $V$ & $S$ & $R$ \\ 
No Rejection & $U$ & $T$ & $m-R$ \\ \hline
Total & $m_0$ & $m_1$ & $m$
\end{tabular}
\end{table}

A single test is evaluated based on its probabilities of type I and type II errors, or equivalently, its size and power.  A good test is a test with small size and large power.  Analogous measures need to be considered when evaluating a multiple test: one based on the number of rejected true nulls and another one based on the number of false nulls rejected ($V$ and $S$ in \tabref{multiple}).  A good multiple test is a multiple test with small number of rejected true nulls and large number of rejected false nulls.    
The mFDR provides a notion of size for multiple tests.  Given a multiple test (left implicit) resulting in the situation of \tabref{multiple}, it is defined as 
\beq
\mfdr = \frac{\E(V)}{\E(R)}.
\eeq

\citet{ioannidis2005most} implicitly considers the $\alpha$-threshold multiple test, defined by
\beq\label{multiple-alpha}
\cR_\alpha = \{i : p_i \le \alpha\},
\eeq
where $\alpha$ is given.  

\begin{prp}
Assume that each test has power $1-\beta$ at level $\alpha$.
Then  
\beq
\ppv = 1 - \mfdr,
\eeq
where $\mfdr$ is the mFDR for the $\alpha$-threshold multiple test \eqref{multiple-alpha}.
\end{prp}

The proof is straightforward.
Indeed, the expected number of rejections is 
\beq
\E(R) = m_0 \alpha + m_1 (1-\beta),
\eeq
and the expected number of incorrect rejections is 
\beq
\E(V) = m_0 \alpha.
\eeq
Thus 
\beq
\mfdr = \frac{m_0 \alpha}{m_0 \alpha + m_1 (1-\beta)} = 1 - \frac{\pi_1 (1-\beta)}{\pi_0 \alpha + \pi_1 (1-\beta)},
\eeq
after dividing by $m$ in the numerator and denominator, and we conclude with the fact that $\pi_0 = 1 - \pi_1$.

\section{Optimal multiple testing procedures}
The last section gives a multiple testing interpretation of the PPV as one minus the mFDR that results from using the $\alpha$-threshold multiple test \eqref{multiple-alpha} in a context where the tests have power $1-\beta$ at level $\alpha$.  
We already discussed in \eqref{ppv-ub} some implications for using this multiple test, but now that we placed the problem in the larger context of multiple testing, one has to wonder if a more sophisticated approach would improve the situation.
Indeed, in our utopian setting, the field as a whole has all interest in using a good multiple test to decide, in a principled way, which studies are significant and worth of further investigation.  

Surely, there is a large and growing statistical literature on the topic of multiple testing, offering a number of methods.  We refer the reader to some books~\cite{dudoit2007multiple,dickhaus2014simultaneous} and some papers surveying the field~\cite{roquain2011type,finner2002partitioning}.
%The literature is extensive in large part because it considers more complicated settings where P-values can exhibit various kinds of dependency.  The setting here is the simplest possible, as we assume that the P-values are independent.  

\subsection{The family-wise error rate (FWER)}
The family-wise error rate (FWER) of a given multiple test is defined as the probability of an incorrect rejection, or in the notation of \tabref{multiple},
\beq
\fwer = \P(V \ge 1).
\eeq

The $\alpha$-threshold multiple test is still popular, for example, in genetics.  It is discussed in~\cite{barsh2012guidelines} in the context of genome-wide association studies (GWAS), where $\alpha$ is said to be typically set to $5 \times 10^{-8}$.  The rationale is that a typical study probes on the order of $10^6$ (one million) SNP's, so that the threshold is set by Bonferroni's method ($0.05/10^6 = 5 \times 10^{-8}$).  This guaranties control of the FWER in a given study.
However, it does not control for the fact that hundreds, if not thousands of such studies are performed each year.

\subsection{The false discovery rate (FDR)}
For many decades the FWER was the main measure of size for multiple tests and apparently remains quite popular.  However, in contexts where a large number of tests are performed (large $m$), controlling the FWER is known to severely limit the (multiple) power.  
This lead~\citet{benjamini1995controlling} to introduce the false discovery rate (FDR), which in the notation of \tabref{multiple} is defined as
\beq
\fdr = \E\left(\frac{V}{\max(R, 1)}\right).
\eeq
This is now a popular notion of size for multiple tests, for example, in genetics~\cite{dudoit2007multiple,bush2012genome}.  

The mFDR is a later variant of the FDR~\cite{storey2007optimal}.
Although the mFDR is conceptually simpler than the FDR, the former cannot be properly controlled when all the null hypotheses are true ($m_0 = m$) while it is still possible to control the latter.
That said, they are known to be equivalent under somewhat mild assumptions when the number of tests is large.
Nevertheless, the FDR appears more popular in practice and we consider control of the FDR in what follows.

\subsection{The Benjamini-Hochberg (BH) procedure in the normal means model}
\tabref{BH} describes the famous procedure of~\citet{benjamini1995controlling}, known to control the FDR (at a prescribed $\gamma$) when the P-values are independent.\footnote{The method is also known to control the FDR at the prescribed level under some forms of positive dependence.}
Note that this multiple test is of threshold type, meaning of the form
\beq\label{threshold}
\cR = \{i : p_i \le \tau\}.
\eeq
For the BH multiple test, the threshold is 
\beq\label{taubh}
\taubh = \hat s \gamma/m, \quad \text{where $\hat s$ is defined in \tabref{BH}.}
\eeq
Note that it is a function of the P-values.

\begin{table}
\centering
\caption{The~\citet{benjamini1995controlling} (BH) multiple test.}
\label{tab:BH}
\bigskip
\begin{tabular}{p{0.7\textwidth}}
\toprule
$p_1, \dots, p_m$ denote the P-values;
$\gamma$ denotes the desired FDR level \\
\midrule
{\bf 1:} Order the P-values, obtaining $p_{(1)} \le \cdots \le p_{(m)}$. \\
{\bf 2:} Let $\hat s$ denote the largest $s$ such that $p_{(s)} \le s \gamma/m$. \\ 
{\bf 3:} Reject the nulls corresponding to the $\hat s$ smallest P-values. \\
\bottomrule
\end{tabular}
\end{table}

%This method is now well-understood.  In particular, the following enlightening result is available.
%\begin{prp}[\citet{genovese2002operating}]
%Suppose that the P-values are independent, and that they are identically distributed under the null and alternative.  Specifically, we assume that $p_i$ is uniformly distributed in $[0,1]$ under its null $\H_i$ and has distribution function $F$ under its alternative $\H_i^\comp$.  Define 
%\beq
%\lambda = \frac{1-\alpha + \alpha \pi_1}{\alpha \pi_1},
%\eeq
%and assume that there is a unique point $u_* \in [0,1]$ such that $F(u_*) = \lambda u_*$ and $F'(u_*) \ne \beta$.  Then in an asymptotic setting where $m\to \infty$ while $\pi_1$ tends to a positive limit (still denoted $\pi_1$), and everything else remaining fixed, the threshold $\tau$ converges to $u_*$.
%\end{prp}

%\subsection{The normal means model}
The normal means model assumes that each test being performed is a Z-test and that any underlying effect only affects the mean.  In detail, we assume that for each $i \in \{1, \dots, m\}$, we have available a Z-ratio $X_i$ and that we reject for large values of $X_i$.  The P-value is computed assuming that $X_i$ is standard normal under $\H_i$, meaning $p_i = 1 - \Phi(X_i)$, where $\Phi$ denotes the standard normal distribution function.  We further assume that the model is indeed accurate in that $X_i \sim \cN(0,1)$ under $\H_i$, while $X_i \sim \cN(\mu, 1)$ under $\H_i^\comp$, where $\mu > 0$ (since we are rejecting for large values of $X_i$).
For example,~\citet{ingster1997some} and~\citet{donoho2004higher} considered the problem of testing the global null hypothesis ($\bigcap_{i=1}^m \H_i$) in this setting.
Note that in this model $p_i \sim \Unif(0,1)$ under $\H_i$, while $p_i \sim F_\mu$ under $\H_i^\comp$, where $F_\mu(u) := 1 - \Phi(\Phi^{-1}(1-u) - \mu)$.

In an asymptotic setting where $m \to \infty$ but everything else remains essentially constant, we have the following consequence of a more general result of~\citet{genovese2002operating}.
We first remark that for any $\mu > 0$ and any $\lambda > 1$, there is a unique point $u \in (0,1)$ such that $F_\mu(u) = \lambda u$, which we denote by $u(\mu, \lambda)$.  We note that $u(\mu, \lambda)$ is a decreasing function of $\lambda$ and that $u(\mu, \infty) = 0$.

\begin{prp}\label{prp:GW} %[\citet{genovese2002operating}] 
Consider an asymptotic setting where $m\to \infty$ while $\pi_1 = \pi_1(m) \to \pi_1^\infty \in (0,1)$ and $\mu = \mu(m) \to \mu^\infty > 0$.  Then $\taubh$ defined in \eqref{taubh} converges in probability to $u(\mu^\infty, \lambda^\infty)$ where $\lambda^\infty = 1 + (1-\gamma)/\gamma \pi_1^\infty$.
\end{prp}

This result implies that, as the proportion of false nulls $\pi_1$ becomes negligible ($\pi_1 \to 0$), the threshold for the BH multiple test converges to zero ($\tau_\gamma \to 0$) .  The conclusion is thus morally the same as with the $\alpha$-threshold multiple test: a small $\pi_1$ makes rejection more difficult, that is, if we use the BH procedure for controlling the FDR.

\subsection{Optimality in the normal means model}
It is possible to study questions of optimality in the context of multiple testing in the normal means model.  
For example, \citet{jin2014rare} do this in a more general minimax framework\footnote{The setting in that paper is that of classification, which is intimately related to the problem of multiple testing.} that includes the present context. 
We present here an alternative, and arguably simpler approach, that based on oracle information~\cite{sun2007oracle,meinshausen2011asymptotic,arias2016distribution}.

%In recent work~\cite{arias2016distribution} we derive an oracle bound instead.  
Suppose that an oracle tells us exactly which null hypotheses are false.  With that information, we would just reject precisely these hypotheses.  But knowing this does not enlighten the discussion!  Imagine, however, that we are restricted to use a threshold procedure, meaning that the rejection set needs to be of the form \eqref{threshold}.  This means that if we reject $\H_i$ then we must reject $\H_j$ for all $j$ such that $p_j \le p_i$.  Remember that in our setting there is nothing that distinguishes the null hypotheses a priori, so doing this is rather natural.  
%In the present context, any reasonable procedure is of the form 
%\beq\label{T}
%\text{declare $p_j$ significant if $p_j \le \tau$},
%\eeq
%where $\tau = \tau(p_1, \dots, p_m)$ is a threshold that may be a function of the P-values.  
The $\alpha$-threshold multiple test (implicitly considered by Ioannidis) and the Benjamini-Hochberg procedure are certainly of this form.  
Now, what is the best that we can do with this oracle information (i.e., full disclosure) under the constraint that we must use a threshold procedure?
A good candidate is the largest threshold that controls the false discovery proportion (FDP) at the desired level.  For a given multiple test, in the notation of \tabref{multiple},
\beq
\fdp = \frac{V}{\max(R, 1)},
\eeq
so that $\fdr = \E(\fdp)$.  Let $V_t$ and $R_t$, and $\fdp_t$, denote the corresponding quantities when we select the P-values not exceeding  $t$.  Define the oracle threshold procedure as
\beq
\tau_\gamma^* = \max\big\{t : \fdp_t \le \gamma\big\}.
\eeq
(We note that the procedure we consider in~\cite{arias2016distribution} is slightly different.)

Following the same arguments leading to \prpref{GW}, we can derive the asymptotic behavior of the oracle threshold multiple test.

\begin{prp}
In the same setting as \prpref{GW}, the oracle threshold $\tau_\gamma^*$ also converges in probability to the same limit $u(\mu^\infty, \lambda^\infty)$.
\end{prp} 

Thus the BH procedure achieves the same first-order asymptotic performance as the oracle procedure.  This is true in a setting where the model parameters remain constant as the number of hypotheses increases.  
In~\cite{arias2016distribution} we consider the situation where the model parameters are allowed to change with $m$ and arrive at essentially the same conclusion.\footnote{The BH procedure is likely optimal also in a minimax sense, although this has not been proved.}

If the oracle procedure is the best we can hope to achieve, then we can conclude, as before, that in order to control the FDR the threshold for rejection needs to be taken closer and closer to 0 as the number of hypotheses increases.  Because of that, there is a real danger that the false nulls drown in the midst of a large number of true null hypotheses.

\section{Discussion} 
\label{sec:discussion}
By properly placing the PPV in the context of multiple testing, we provided some additional theoretical backing to Ioannidis' claim \eqref{ioannidis}.  Perhaps it is no longer just about what each lab does in its corner; perhaps it also matters what other labs are doing.  For if some labs are repeatedly probing the world more or less at random, without much thought, the hypotheses they test and the P-values they publish will act as noise and possibly overwhelm the field to a point where it becomes impossible, at least with statistical methodology alone, to separate the false nulls from the true nulls.

\medskip\noindent {\em Field-wide control of false discovery rate?}
At a higher level, we question whether research can continue to be done independently by a large number of researchers without some sort of central decision making on what constitutes important questions.  
There are some initiatives in this direction~\cite{open2012open}.  
\citet{rosset2014novel}, very recently, proposed a control via access to a large database.  Note that this proposal requires the collection of new samples so that the database is continually enriched and FDR control remains possible.

\medskip\noindent {\em Blame it on the P-values?}
With the bad reputation that P-values are getting and the long-standing controversy around them~\cite{nuzzo2014statistical,wasserstein2016asa,murtaugh2014defense,trafimow2015editorial}, one may wonder if all this can also be blamed on them.  For example, one may wonder if the use of Bayes factors~\cite{goodman1999toward-part1,goodman1999toward-part2}, or even confidence intervals, would solve the problem.  The answer is a resounding {\em no} if, as we consider here, the goal is to identify the false nulls with some degree of accuracy.

\subsection*{Acknowledgments}
I am grateful to my colleagues, Jeff Rabin and Ian Abramson, for  extensive discussions on the topic and for reading a draft of the manuscript.  
This work was partially supported by a grant from the US National Science Foundation (DMS 1223137).

\small
\bibliographystyle{plainnat}
\bibliography{observer-effect}

\end{document}